\documentstyle[11pt]{article}

\begin{document}

\noindent{\bf NONCLASSICAL TYPE REPRESENTATIONS OF THE $q$-DEFORMED
ALGEBRA $U'_q({\rm so}_n)$}
\bigskip

\centerline{N. Z. Iorgov, A. U. Klimyk}
\medskip

\centerline{Institute for Theoretical Physics,
Kiev 252143, Ukraine}
\bigskip

\begin{abstract}
The nonstandard $q$-deformation $U'_q({\rm so}_n)$ of the universal
enveloping algebra $U({\rm so}_n)$ has irreducible finite dimensional
representations which are a $q$-deformation of the well-known
irreducible finite dimensional representations of $U({\rm so}_n)$.
But $U'_q({\rm so}_n)$ also has irreducible finite dimensional
representations which have no classical analogue. The aim of this
paper is to give these representations which are called nonclassical
type representations. They are given by explicit formulas for
operators of the representations corresponding to the generators of
$U'_q({\rm so}_n)$.

\end{abstract}

{\bf 1.} In [1] it was constructed a $q$-deformation $U'_q({\rm so}_n)$
of the universal enveloping algebra $U({\rm so}_n)$ which differ from
the quantum algebra $U_q({\rm so}_n)$ introduced by V. Drinfeld [2]
and M. Jimbo [3] (see also [4]). The algebra $U'_q({\rm so}_n)$
permits to construct the reduction of $U'_q({\rm so}_n)$ onto
$U'_q({\rm so}_{n-1})$ which can be used for construction of an
analogue of Gel'fand--Tsetlin bases for irreducible representations.

In the classical case, the imbedding $SO(n)\subset SU(n)$
(and its infinitesimal analogue) is of great importance for nuclear
physics and in the theory of Riemannian spaces. It is well known
that in the framework of Drinfeld--Jimbo quantum groups and algebras
one cannot construct the corresponding imbedding. The algebra
$U'_q({\rm so}_n)$ allows to define such an imbedding [5],
that is, it is possible to define the imbedding
$U'_q({\rm so}_n)\subset U_q({\rm sl}_n)$,
where $U_q({\rm sl}_n)$ is the Drinfeld-Jimbo quantum algebra.

As a disadvantage of the algebra $U'_q({\rm so}_n)$ we have
to mention the difficulties with Hopf algebra structure. Nevertheless,
$U'_q({\rm so}_n)$ turns out to be a coideal in
$U_q({\rm sl}_n)$.

Finite dimensional irreducible representations of the algebra
$U'_q({\rm so}_n)$ were constructed in [1]. The formulas
of action of the generators of $U'_q({\rm so}_n)$ upon the
basis (which is a $q$-analogue of the Gel'fand--Tsetlin basis) are
given there. A proof of these formulas and some their corrections were
given in [6]. However,
finite dimensional irreducible representations described in [1] and [6]
are representations of the classical type. They are $q$-deformations of the
corresponding irreducible representations of the Lie algebra
${\rm so}_n$, that is, at $q\to 1$ they turn into representations
of ${\rm so}_n$.

The algebra $U'_q({\rm so}_n)$ has other classes of finite
dimensional irreducible representations which have no classical analogue.
These representations are singular at the limit $q\to 1$. The aim
of this paper is to describe these representations of
$U'_q({\rm so}_n)$. Note that the description of these
representations for the algebra $U'_q({\rm so}_3)$ is given in
[7]. A classification of irreducible $*$-representations of real forms
of the algebra $U'_q({\rm so}_3)$ is given in [8].

The algebra $U'_q({\rm so}_n)$
is defined as the associative algebra (with a unit)
generated by the elements $I_{i,i-1}$, $i=2,3,...,n$,
satisfying the defining relations
$$
I_{i+1,i}I^2_{i,i-1}-(q+q^{-1})I_{i,i-1}I_{i+1,i}I_{i,i-1}
+I^2_{i,i-1}I_{i+1,i}=-I_{i+1,i},     \eqno (1)
$$
$$
I^2_{i+1,i}I_{i,i-1}-(q+q^{-1})I_{i+1,i}I_{i,i-1}
I_{i+1,i}+I_{i,i-1}I^2_{i+1,i}=-I_{i,i-1},   \eqno(2)
$$
$$
[I_{i,i-1},I_{j,j-1}]=0,\ \ \ \vert i-j\vert >1,    \eqno (3)
$$
where [.,.] denotes the usual commutator.
In the limit $q\to 1$ formulas (1)--(3) give the relations
defining the universal enveloping algebra $U({\rm so}_n)$.

Note also that relations
(1) and (2) principally differ from the $q$-deformed Serre relations
in the approach of Drinfeld [2] and Jimbo [3] to quantum orthogonal
algebras by a presence of nonzero right hand side and by possibility
of the reduction
$U'_q({\rm so}_n)\supset U'_q({\rm so}_{n-1})\supset
\cdots \supset U'_q({\rm so}_3)$.
\medskip

{\bf 2.} In this section we describe
irreducible finite dimensional
representations of the algebras $U'_q({\rm so}_{n})$, $n \ge 3$,
which are $q$-deformations of the finite dimensional irreducible
representations of the Lie algebra ${\rm so}_n$.
They are given by sets ${\bf m}_{n}$
consisting of $\{ {n/2}\}$ numbers $m_{1,n}, m_{2,n},..., m_{\left
\{ {n/2}\right \} ,n}$ (here $\{ {n/2}\}$ denotes integral part
of ${n/2}$) which are all integral or all half-integral and
satisfy the dominance conditions
$$
m_{1,2p+1}\ge m_{2,2p+1}\ge ... \ge
m_{p,2p+1}\ge 0 ,                                  \eqno(4)
$$
$$
m_{1,2p}\ge m_{2,2p}\ge ... \ge m_{p-1,2p}\ge |m_{p,2p}|
                                                    \eqno(5)
$$
for $n=2p+1$ and $n=2p$, respectively.
These representations are denoted by $T_{{\bf m}_n}$.
For a basis in a representation space we take the $q$-analogue of the
Gel'fand--Tsetlin basis which is obtained by successive reduction of
the representation $T_{{\bf m}_n}$ to the subalgebras
$U'_q({\rm so}_{n-1})$, $U'_q({\rm so}_{n-2})$, $\cdots$, $U'_q({\rm
so}_3)$,
$U'_q({\rm so}_2):=U({\rm so}_2)$.
As in the classical case, its elements are labelled by Gel'fand--Tsetlin
tableaux
$$
  \{\xi_{n} \}
\equiv \left\{ \matrix{ {\bf m}_{n} \cr {\bf m}_{n-1} \cr \dots
\cr {\bf m}_{2}  }
 \right\}
\equiv \{ {\bf m}_{n},\xi_{n-1}\}\equiv \{{\bf m}_{n} ,
{\bf m}_{n-1} ,\xi_{n-2}\} ,
                                                       \eqno(6)
$$
where the components of ${\bf m}_{n}$ and ${\bf m}_{n-1}$ satisfy the
"betweenness" conditions
$$
m_{1,2p+1}\ge m_{1,2p}\ge m_{2,2p+1} \ge m_{2,2p} \ge ...
\ge m_{p,2p+1} \ge m_{p,2p} \ge -m_{p,2p+1}  ,
                                                     \eqno(7)
$$
$$
m_{1,2p}\ge m_{1,2p-1}\ge m_{2,2p} \ge m_{2,2p-1} \ge ...
\ge m_{p-1,2p-1} \ge \vert m_{p,2p} \vert .
                                                      \eqno(8)
$$
The basis element defined by tableau $\{\xi_{n} \}$ is denoted
as $\vert \{\xi_{n} \} \rangle $ or simply as $\vert
\xi_{n} \rangle $.

It is convenient to introduce the so-called $l$-coordinates
$$
l_{j,2p+1}=m_{j,2p+1}+p-j+1,  \qquad
                      l_{j,2p}=m_{j,2p}+p-j ,        \eqno(9)
$$
for the numbers $m_{i,k}$.
In particular, $l_{1,3}=m_{1,3}+1$ and $l_{1,2}=m_{1,2}$.
The operator $T_{{\bf m}_n}(I_{2p+1,2p})$ of the representation
$T_{{\bf m}_n}$ of $U'_q({\rm so}_{n})$ acts upon Gel'fand--Tsetlin
basis elements, labelled by (6), by the formula
$$
T_{{\bf m}_n}(I_{2p+1,2p})
| \xi_n\rangle =
\sum^p_{j=1} \frac{ A^j_{2p}(\xi_n)}
{q^{l_{j,2p}}+q^{-l_{j,2p}} }
            \vert (\xi_n)^{+j}_{2p}\rangle -
\sum^p_{j=1} \frac{A^j_{2p}((\xi_n)^{-j}_{2p})}
{q^{l_{j,2p}}+q^{-l_{j,2p}}}
|(\xi_n)^{-j}_{2p}\rangle                         \eqno(10)
$$
and the operator $T_{{\bf m}_n}(I_{2p,2p-1})$ of the representation
$T_{{\bf m}_n}$ acts as
$$
T_{{\bf m}_n}(I_{2p,2p-1})\vert \xi_n\rangle=
\sum^{p-1}_{j=1} \frac{B^j_{2p-1}(\xi_n)}
{[2 l_{j,2p-1}-1][l_{j,2p-1}]}
\vert (\xi_n)^{+j}_{2p-1} \rangle
$$
$$
-\sum^{p-1}_{j=1}\frac {B^j_{2p-1}((\xi_n)^{-j}_{2p-1})}
{[2 l_{j,2p-1}-1][l_{j,2p-1}-1]}
\vert (\xi_n)^{-j}_{2p-1}\rangle
+ {\rm i}\, C_{2p-1}(\xi_n)
\vert \xi_n \rangle .                                 \eqno(11)
$$
In these formulas, $(\xi_n)^{\pm j}_{k}$ means the tableau (6)
in which $j$-th component $m_{j,k}$ in ${\bf m}_k$ is replaced
by $m_{j,k}\pm 1$. The coefficients
$A^j_{2p},  $ $B^j_{2p-1},$ $C_{2p-1}$ in (10) and (11) are given
by the expressions
$$
\leqno A^j_{2p}(\xi_n) =
$$
$$
\left( \frac{\prod_{i=1}^p [l_{i,2p+1}+l_{j,2p}] [l_{i,2p+1}-l_{j,2p}-1]
\prod_{i=1}^{p-1} [l_{i,2p-1}+l_{j,2p}] [l_{i,2p-1}-l_{j,2p}-1]}
{\prod_{i\ne j}^p [l_{i,2p}+l_{j,2p}][l_{i,2p}-l_{j,2p}]
[l_{i,2p}+l_{j,2p}+1][l_{i,2p}-l_{j,2p}-1]} \right)^{1\over 2}  \eqno (12)
$$
and
$$
\leqno  B^j_{2p-1}(\xi_n)=
$$
$$
\left( \frac{\prod_{i=1}^p
[l_{i,2p}+l_{j,2p-1}] [l_{i,2p}-l_{j,2p-1}] \prod_{i=1}^{p-1}
[l_{i,2p-2}+l_{j,2p-1}] [l_{i,2p-2}-l_{j,2p-1}]}
{\prod_{i\ne j}^{p-1}
[l_{i,2p-1}{+}l_{j,2p-1}][l_{i,2p-1}{-}l_{j,2p{-}1}]
[l_{i,2p-1}{+}l_{j,2p-1}{-}1][l_{i,2p-1}{-}l_{j,2p-1}{-}1]}
 \right) ^{1\over 2} ,
\eqno(13)
$$
$$
C_{2p-1}(\xi_n) =\frac{ \prod_{s=1}^p [ l_{s,2p} ]
\prod_{s=1}^{p-1} [ l_{s,2p-2} ]}
{\prod_{s=1}^{p-1} [l_{s,2p-1}] [l_{s,2p-1} - 1] } ,   \eqno(14)
$$
where numbers in square brackets mean $q$-numbers defined by
$$
[a]:= (q^a-q^{-a})/(q-q^{-1}).
$$
It is seen from (9) that $C_{2p-1}$ in (14) identically
vanishes if $m_{p,2p}\equiv l_{p,2p}=0$.

A proof of the fact that formulas (10)--(14) indeed determine
a representation of $U'_q({\rm so}_n)$ is given in [6].
\medskip

{\bf 3.}
The representations of the previous section are called representations
of the classical type, since under the limit $q\to 1$ the operators
$T_{{\bf m}_n}(I_{j,j-1})$ turn into the corresponding operators
$T_{{\bf m}_n}(I_{j,j-1})$ for irreducible finite dimensional
representations with highest weights ${\bf m}_n$ of the Lie algebra
${\rm so}_n$.

The algebra $U'_q({\rm so}_n)$ also has irreducible finite dimensional
representations $T$ of nonclassical type, that is, such that the operators
$T(I_{j,j-1})$ have no classical limit $q\to 1$.
They are given by sets $\epsilon := (\epsilon _2,\epsilon _3,\cdots ,
\epsilon _n)$, $\epsilon _i=\pm 1$, and by sets
${\bf m}_{n}$ consisting of $\{ {n/2}\}$ {\bf half-integral} numbers
$m_{1,n}, m_{2,n}, \cdots , m_{\{ n/2\} ,n}$
(here $\{ {n/2}\}$ denotes integral part of ${n/2}$) that satisfy the
dominance conditions
$$
m_{1,2p+1}\ge m_{2,2p+1}\ge ... \ge
m_{p,2p+1}\ge 1/2 ,                                  \eqno(15)
$$
$$
m_{1,2p}\ge m_{2,2p}\ge ... \ge m_{p-1,2p}\ge m_{p,2p}\ge 1/2
                                                    \eqno(16)
$$
for $n=2p+1$ and $n=2p$, respectively.
These representations are denoted by $T_{\epsilon,{\bf m}_n}$.

For a basis in the representation space we use the analogue of the
basis of the previous section. Its elements are
labelled by tableaux
$$
  \{\xi_{n} \}
\equiv \left\{ \matrix{ {\bf m}_{n} \cr {\bf m}_{n-1} \cr \dots
\cr {\bf m}_{2}  }
 \right\}
\equiv \{ {\bf m}_{n},\xi_{n-1}\}\equiv \{{\bf m}_{n} ,
{\bf m}_{n-1} ,\xi_{n-2}\} ,
                                                       \eqno(17)
$$
where the components of ${\bf m}_{2p+1}$ and ${\bf m}_{2p}$ satisfy the
"betweenness" conditions
$$
m_{1,2p+1}\ge m_{1,2p}\ge m_{2,2p+1} \ge m_{2,2p} \ge ...
\ge m_{p,2p+1} \ge m_{p,2p} \ge 1/2  ,
                                                     \eqno(18)
$$
$$
m_{1,2p}\ge m_{1,2p-1}\ge m_{2,2p} \ge m_{2,2p-1} \ge ...
\ge m_{p-1,2p-1} \ge m_{p,2p}  .
                                                      \eqno(19)
$$
The basis element defined by tableau $\{\xi_{n} \}$ is denoted
as $\vert \xi_{n} \rangle $.

As in the previous section, it is convenient to introduce the
$l$-coordinates
$$
l_{j,2p+1}=m_{j,2p+1}+p-j+1,  \qquad
                      l_{j,2p}=m_{j,2p}+p-j .        \eqno(20)
$$
The operator $T_{\epsilon,{\bf m}_n}(I_{2p+1,2p})$ of the representation
$T_{\epsilon,{\bf m}_n}$ of $U_q({\rm so}_{n})$ acts upon our
basis elements, labelled by (17), as
$$
T_{\epsilon,{\bf m}_n}(I_{2p+1,2p})
| \xi_n\rangle =
\delta_{m_{p,2p},1/2}\, \frac{\epsilon _{2p+1}}{q^{1/2}-q^{-1/2}} D_{2p}
(\xi _n) | \xi_n\rangle +
$$
$$
+\sum^p_{j=1} \frac{ A^j_{2p}(\xi_n)}
{q^{l_{j,2p}}-q^{-l_{j,2p}} }
            \vert (\xi_n)^{+j}_{2p}\rangle -
\sum^p_{j=1} \frac{A^j_{2p}((\xi_n)^{-j}_{2p})}
{q^{l_{j,2p}}-q^{-l_{j,2p}}}
|(\xi_n)^{-j}_{2p}\rangle    ,                     \eqno(21)
$$
where $\delta_{m_{p,2p},1/2}$ is the Kronecker symbol and
the summation in the last sum must be from 1 to $p-1$ if
$m_{p,2p}=1/2$.
The operator $T_{\epsilon ,{\bf m}_n}(I_{2p,2p-1})$ of the representation
$T_{\epsilon ,{\bf m}_n}$ acts as
$$
T_{\epsilon ,{\bf m}_n}(I_{2p,2p-1})\vert \xi_n\rangle=
\sum^{p-1}_{j=1} \frac{B^j_{2p-1}(\xi_n)}
{[2 l_{j,2p-1}-1][l_{j,2p-1}]_+}
\vert (\xi_n)^{+j}_{2p-1} \rangle -
$$
$$
-\sum^{p-1}_{j=1}\frac {B^j_{2p-1}((\xi_n)^{-j}_{2p-1})}
{[2 l_{j,2p-1}-1][l_{j,2p-1}-1]_+}
\vert (\xi_n)^{-j}_{2p-1}\rangle
+ \epsilon _{2p} {\hat C}_{2p-1}(\xi_n)
\vert \xi_n \rangle ,                                 \eqno(22)
$$
where
$$
[a]_+=(q^a+q^{-a})/(q-q^{-1}).
$$
In these formulas, $(\xi_n)^{\pm j}_{k}$ means the tableau (17)
in which $j$-th component $m_{j,k}$ in ${\bf m}_k$ is replaced
by $m_{j,k}\pm 1.$ The coefficients $A^j_{2p}$ and $B^j_{2p-1}$
in (21) and (22) are given by the same formulas
as in (10) and (11) (that is, by the formulas (12) and (13)) and
$$
{\hat C}_{2p-1}(\xi_n) = {
\prod_{s=1}^p [ l_{s,2p} ]_+
\prod_{s=1}^{p-1} [ l_{s,2p-2} ]_+ \over
   \prod_{s=1}^{p-1} [l_{s,2p-1}]_+ [l_{s,2p-1} - 1]_+ } ,   \eqno(23)
$$
$$
D_{2p} (\xi _n)=
\frac{\prod_{i=1}^p
[l_{i,2p+1}-\frac 12 ] \prod_{i=1}^{p-1} [l_{i,2p-1}-\frac 12 ] }
{\prod_{i=1}^{p-1}
[l_{i,2p}+\frac 12 ] [l_{i,2p}-\frac 12 ] } . \eqno (24)
$$

The fact that the above operators $T_{\epsilon ,{\bf m}_n}(I_{k,k-1})$
satisfy the defining relations (1)--(3) of the algebra $U'_q({\rm so}_n)$
is proved in the following way. We take the formulas (10)--(14) for
the classical type representations $T_{{\bf m}_n}$ of $U'_q({\rm so}_n)$
with half-integral $m_{i,n}$ and replace there every
$m_{j,2p+1}$ by $m_{j,2p+1}-{\rm i}\pi /2h$,
every $m_{j,2p}$, $j\ne p$, by $m_{j,2p}-{\rm i}\pi /2h$
and $m_{p,2p}$ by $m_{p,2p}-\epsilon _2 \epsilon _4\cdots \epsilon _{2p}
{\rm i}\pi /2h$, where each $\epsilon _{2s}$ is equal
to $+1$ or $-1$ and $h$ is defined by $q=e^h$.
Repeating almost word by word the reasoning of the paper [6],
we prove that the operators given by
formulas (10)--(14) satisfy the defining relations (1)--(3)
of the algebra $U'_q({\rm so}_n)$ after this replacement. Therefore,
these operators determine a representation of $U'_q({\rm so}_n)$. We
denote this representation by $T'_{{\bf m}_n}$.
After a simple rescaling, the operators
$T'_{{\bf m}_n}(I_{k,k-1})$ take the form
$$
T'_{{\bf m}_n}(I_{2p+1,2p})
| \xi_n\rangle =
\sum^p_{j=1} \frac{ A^j_{2p}(\xi_n)}
{q^{l_{j,2p}}-q^{-l_{j,2p}} }
            \vert (\xi_n)^{+j}_{2p}\rangle -
\sum^p_{j=1} \frac{A^j_{2p}((\xi_n)^{-j}_{2p})}
{q^{l_{j,2p}}-q^{-l_{j,2p}}}
|(\xi_n)^{-j}_{2p}\rangle    ,
$$
$$
T'_{{\bf m}_n}(I_{2p,2p-1})\vert \xi_n\rangle=
\sum^{p-1}_{j=1} \frac{B^j_{2p-1}(\xi_n)}
{[2 l_{j,2p-1}-1][l_{j,2p-1}]_+}
\vert (\xi_n)^{+j}_{2p-1} \rangle -
$$
$$
-\sum^{p-1}_{j=1}\frac {B^j_{2p-1}((\xi_n)^{-j}_{2p-1})}
{[2 l_{j,2p-1}-1][l_{j,2p-1}-1]_+}
\vert (\xi_n)^{-j}_{2p-1}\rangle
+ \epsilon _{2p} {\hat C}_{2p-1}(\xi_n)
\vert \xi_n \rangle ,
$$
where $A^j_{2p}$, $B^j_{2p-1}$ and ${\hat C}_{2p-1}$ are such as in
the formulas (21) and (22).
The representations
$T'_{{\bf m}_n}$ are reducible. We decompose these representations into
subrepresentations in the following way.
We fix $p$ ($p=1,2,\cdots ,\{ (n-1)/2\}$)
and decompose the linear space ${\cal H}$ of the
representation $T'_{{\bf m}_n}$ into direct sum of two
subspaces ${\cal H}_{\epsilon _{2p+1}}$,
$\epsilon _{2p+1}=\pm 1$, spanned by the basis vectors
$$
| \xi_n \rangle _{\epsilon _{2p+1}} =| \xi_n \rangle -\epsilon _{2p+1}
| \xi_n' \rangle , \ \ \ \ \ m_{p,2p}\ge 1/2,
$$
respectively,
where $| \xi_n' \rangle$ is obtained from $| \xi_n \rangle$ by replacement
of $m_{p,2p}$ by $-m_{p,2p}$.
A direct verification shows that two subspaces
${\cal H}_{\epsilon _{2p+1}}$ are
invariant with respect to all the operators $T'_{{\bf m}_n}(I_{k,k-1})$.
Now we take the subspaces ${\cal H}_{\epsilon _{2p+1}}$ and repeat the
same procedure for some $s$, $s\ne p$, and decompose each of these
subspaces into two invariant subspaces. Continuing this procedure further
we decompose the representation space ${\cal H}$ into a direct sum
of $2^{\{ (n-1)/2\} }$ invariant subspaces. The operators
$T'_{{\bf m}_n}(I_{k,k-1})$ act upon these subspaces by the formulas
(21) and (22). We denote the corresponding subrepresentations on these
subspaces by $T_{\epsilon ,{\bf m}_n}$. The above reasoning shows that
the operators
$T_{\epsilon ,{\bf m}_n}(I_{k,k-1})$ satisfy the defining relations
(1)--(3) of the algebra $U'_q({\rm so}_n)$.
\medskip

\noindent
{\bf Theorem 1.} {\it The representations $T_{\epsilon ,{\bf m}_n}$
are irreducible.
The representations $T_{\epsilon ,{\bf m}_n}$ and
$T_{\epsilon ',{\bf m}'_n}$ are pairwise nonequivalent for
$(\epsilon ,{\bf m}_n)\ne (\epsilon ',{\bf m}'_n)$.
For any admissible $(\epsilon ,{\bf m}_n)$ and ${\bf m}'_n$
the representations $T_{\epsilon ,{\bf m}_n}$ and $T_{{\bf m}'_n}$ are
pairwise nonequivalent.}
\medskip

The algebra $U'_q({\rm so}_n)$ has non-trivial one-dimensional
representations. They are special cases of the representations of the
nonclassical type. They are described as follows.

Let $\epsilon := (\epsilon _2,\epsilon _3,\cdots ,
\epsilon _n)$, $\epsilon _i=\pm 1$, and let
${\bf m}_{n}=(m_{1,n}, m_{2,n}, \cdots , m_{ \{ n/2 \} ,n})
=(\frac 12 , \frac 12 ,\cdots ,\frac 12 )$. Then the corresponding
representations $T_{\epsilon ,{\bf m}_n} $ are one-dimensional
and are given by the formulas
$$
T_{\epsilon ,{\bf m}_n}(I_{k+1,k})| \xi_n\rangle =
\frac{\epsilon _{k+1}}{q^{1/2}-q^{-1/2}}| \xi_n\rangle .
$$
Thus, to every $\epsilon := (\epsilon _2,\epsilon _3,\cdots ,
\epsilon _n)$, $\epsilon _i=\pm 1$, there corresponds a one-dimensional
representation of $U'_q({\rm so}_n)$.
\medskip

\noindent
{\bf Conjecture.} {\it If $q$ is not a root of unity, then every
irreducible finite dimensional representation of $U'_q({\rm so}_n)$
is equivalent to one of the representations $T_{{\bf m}_n}$ of the
classical type or to one of the representations
$T_{\epsilon ,{\bf m}_n}$ of the nonclassical type.}
\medskip

{\small The research of this publication was made possible in part
by Award No. UP1--309 of CRDF and by Award No. 1.4/206 of Ukrainian
DFFD.

\end{document}